\documentclass{IEEEtran4PSCC}
\ifCLASSINFOpdf
   \usepackage[pdftex]{graphicx}
\else
   \usepackage[dvips]{graphicx}
\fi
\usepackage[cmex10]{amsmath}
\usepackage[style=ieee,backend=bibtex,url=false,doi=false,isbn=false,date=year,citestyle=numeric-comp,maxbibnames=10,minbibnames=10,maxcitenames=10,mincitenames=10]{biblatex}
\bibliography{Aggregation.bib,Convex_models.bib,alexLit}

\usepackage{booktabs}
\usepackage[ruled]{algorithm2e}
\usepackage{xcolor}
\usepackage{amssymb}
\usepackage{subcaption}
\makeatletter
\let\old@ps@headings\ps@headings
\let\old@ps@IEEEtitlepagestyle\ps@IEEEtitlepagestyle
\def\psccfooter#1{%
    \def\ps@headings{%
        \old@ps@headings%
        \def\@oddfoot{\strut\hfill#1\hfill\strut}%
        \def\@evenfoot{\strut\hfill#1\hfill\strut}%
    }%
    \def\ps@IEEEtitlepagestyle{%
        \old@ps@IEEEtitlepagestyle%
        \def\@oddfoot{\strut\hfill#1\hfill\strut}%
        \def\@evenfoot{\strut\hfill#1\hfill\strut}%
    }%
    \ps@headings%
}
\makeatother

\begin{document}
%
\title{Flexibility aggregation via set projection for distribution grids with multiple interconnections }

\author{\IEEEauthorblockN{1\textsuperscript{st} Maísa Beraldo Bandeira}
\IEEEauthorblockA{\textit{ie3 - Institute of
Energy Systems,} \\ \textit{Energy Efficiency and Energy Economics} \\
\textit{TU Dortmund University}\\
Dortmund, Germany \\
maisa.bandeira@tu-dortmund.de}
\and
\IEEEauthorblockN{2\textsuperscript{nd} Alexander Engelmann}
\IEEEauthorblockA{\textit{logarithmo GmbH \& Co. KG} \\
Dortmund, Germany \\
alexander.engelmann@ieee.org}
\and
\IEEEauthorblockN{3\textsuperscript{rd} Timm Faulwasser}
\IEEEauthorblockA{\textit{Institute of Control Systems} \\
\textit{Hamburg University of Technology}\\
Hamburg, Germany \\
timm.faulwasser@ieee.org}
}


\maketitle

\begin{abstract}
With the increasing number of flexible energy devices in distribution grids, coordination between Transmission System Operators (TSOs) and Distribution System Operators (DSOs) becomes critical for optimal system operation. 
One form of coordination is to solve the overall system operation problem in a hierarchical way, computing Feasible Operational Regions (FORs) for the interconnection between TSO/DSO. 
Most methods for computing FORs rely on the assumption of only one interconnection point between TSO and DSOs, which is often violated in practice. In this work, we propose a method for computing FORs in distribution grids with multiple interconnection points to the transmission grid. We test our method in a grid with two interconnecting points and analyze the properties of the resulting high-dimensional FOR from a power systems perspective.
\end{abstract}



\section{Introduction} \label{sec:Intro}
Renewable energy sources and the electrification of the heating and mobility sectors offer large amounts of flexibility, which is a crucial lever for a cost-effective future system operation.  
As many of the corresponding devices are located in distribution grids, they play an increasingly active role in system operation. This raises the question of how to efficiently coordinate TSOs and DSOs. 
Controlling a large number of flexible devices directly at the TSO level is not desirable for many reasons, such as complexity and confidentiality.

Flexibility aggregation has been proposed for reducing complexity and simplifying information exchange. The approach relies on solving the centralized optimization problem in a hierarchical fashion. 
The first step is to compute Feasible Operational Regions (FORs) at the interconnection points between the Transmission System Operator (TSO) and Distributed System Operators (DSOs). 
A FOR is defined as the set of values for the coupling variables, usually active and reactive power, which can be achieved at the interconnection without violating the physical constraints of the distribution grid~\cite{mayorga2018}. 
The TSO uses the FOR in its own optimization problem, leveraging the flexibility of DERs without causing congestion for the DSO. 

Several methods for computing the FOR, also referred to as flexibility aggregation, have been proposed in the literature~\cite{Sarstedt2021}. 
Most of these approaches rely on sampling the feasible region and subsequently computing the convex hull of the resulting points. Different sampling strategies can be employed—for instance, Monte Carlo methods~\cite{MayorgaGonzalez2018}, which, however, require a large number of samples. 
Alternatively, the samples can be determined more systematically by solving an optimal power flow problem to obtain the extreme operating points of the FOR, i.e., the maximum and minimum active and reactive power exchanges at the interconnection~\cite{silva2018, capitanescu2018}. 
These methods perform well for two-dimensional FORs, i.e. for simple PQ-charts. In~\cite{Sarstedt2020}, swarm optimization is used to compute multiple samples that explicitly account for the voltage dependency.


In practice, the interconnection point between grid operators is not always unique. 
In Germany, for example, it is quite common that the TSO and DSO have at least two interconnection points. There are limited works in the literature that address the problem of computing the FOR for more than one interconnection point, since it is not trivial to extend the sampling methods presented above to that scenario. 
The challenge arises from the mutual dependence between the powers at the different interconnections and the physical distance between the buses to which they are connected in the transmission grid, which influences the voltages. 

The work in~\cite{silva2023} extends the method for a grid with two interconnecting points. 
It computes two independent FORs by sampling them individually. 
The coupling between the two interconnection points is considered by deriving a data-based equivalent grid for the transmission grid, but the intrinsic dependency of the flexibility available at each interconnection is neglected. 
Alternatively, the work in~\cite{stark2024} computes an equivalent grid for the distribution grid, rather than a FOR, communicating to the TSO more information than just the feasible values for the coupling variables. However, the authors highlight the importance of the voltage magnitude and angle dependency at the coupling points.

Furthermore, in~\cite{majumdar2023}, the interdependence between the FORs is highlighted. Building on this,~\cite{stark2023} explicitly considers the relationship between the two FORs. They first compute the FOR for one interconnection using an optimization-based sampling procedure. The resulting feasible area is then discretized into a grid. 
For each grid cell—representing a specific range of power exchange at the first interconnection—the method computes a corresponding FOR for the second interconnection, again via sampling. 
The approach yields a set of pairs of FORs for the two interconnections, which is difficult to incorporate directly into the TSO optimization problem.

Recently, set projection has been used for the computation of the FOR~\cite{tan2024, Bandeira2024}. The method has been extended to multi-dimensional FORs, which include time-dependent elements~\cite{wang2023}. In~\cite{lyu2023} set projection was used for the coordination of transmission grids with multiple connections, using a linearized version of the grid that takes into consideration only active power.

However, despite its practical relevance, and to the best of the authors' knowledge, the computation of high-dimensional FORs for distribution grids with multiple interconnections has not been addressed before. Building on our earlier work~\cite{Bandeira2024}, we extend the set projection approach to this setting, enabling the computation of multidimensional FORs that capture all the couplings between interconnection variables. Using a linearized power flow formulation that preserves voltage and reactive power information, the resulting polyhedral FOR can be readily communicated to the TSO as a compact set of inequality constraints. The approach is demonstrated and analyzed for a distribution grid with two interconnection points.

The paper is organized as follows: first, we define our problem and explain our approach for flexibility aggregation in Section \ref{sec:problem}. In Section \ref{sec:simulation}, we present our case study and results. Lastly, Section \ref{sec:conclusion} summarizes our findings and outlines possible future research directions.

\section{Problem Set-up} \label{sec:problem}
We assume the power system to be tree-structured, with the distribution grids having no connections to one another. We define the multi-voltage level power system as a graph  $G^e = (\mathcal N, \mathcal B)$, where each node $ k \in \mathcal N$ is a bus and the edges $(k,l) \in \mathcal B$ represent the branches that connect the buses, e.g., transmission lines or transformers.\\
We consider by $i \in \mathcal S=\{1,\dots,|S|\}$, where $i = 1$ refers to the TSO and the index set 
$\mathcal D \doteq \mathcal S\setminus \{1\}$ refers to the DSOs. 
Furthermore, we decompose the set of buses $\mathcal{N}$ into a set with the TSO buses $\mathcal{N}_1$ and into $|\mathcal D|$  bus sets  for the DSOs $\mathcal N_i, i \in \mathcal{D}$, for all DSOs.
Similarly, 	 we split the set of branches $\mathcal B$ into $\{\mathcal B_i\}_{i \in \mathcal S}$ such that all branches connecting nodes in $\mathcal N_i$ belong to $\mathcal B_i$.
We assign the branches connecting TSO and DSO grids $\mathcal B_i^c \subseteq \mathcal B_i$  to the DSOs.
We denote the buses at the TSO level connected to a branch between TSO and DSO $i$ as $\mathcal N_i^c\subseteq \mathcal N$. Details of the formulation can be found in \cite{Bandeira2025}.

The corresponding OPF problem reads
\begin{subequations}\label{eq:prob1} 
	\begin{align} 
		\min_{\{y_i\}_{i \in \mathcal S},\{z_i\}_{i\in \mathcal D}} &f_1(y_1)+ \sum_{i \in \mathcal D} f_i(y_i,z_i)  \\
		\text{ subject to  } & (y_i,  z_i)\in \mathcal{{ X}}_i \text{ for all }i \in \mathcal D \label{eq:locConstr}\\
		&(y_1,z_2,\dots,z_{|\mathcal S|}) \in \mathcal X_1,
	\end{align}
\end{subequations} where for each subsystem $i$, we denote the local decision variables as $y_i \subseteq x_i$ and the coupling variables as $z_i \subseteq x_i$,\footnote{We use the shorthand  $y\subset x$, $y\in \mathbb R^n,x\in \mathbb R^m, n<m$ to express that $y$ is composed of a subset of the components of $x$, i.e., there exists an index set $\mathcal I \subset \{1,\dots,m\}$ such that $y=[x_i]_{i\in \mathcal I}$, where $[\cdot]$ denotes the vertical  concatenation.} i.e., variables that appear both in the TSO and the DSO problem, so that $x_i = (y_i, z_i)$. Furthermore, $\mathcal X_i$ are local constraint sets of DSOs or the TSO.
 
\subsection{The constraint sets $\mathcal X_i$}
The physical characteristics and limitations of each sub-grid define the feasible sets $\mathcal X_i$ of the subproblems. 
To utilize existing polyhedral set projection tools, we formulate these sets as a convex polyhedron. In \cite{Bandeira2024}, we compare different convex grid models, which lead to polyhedral constraint sets convenient for flexibility aggregation for grids with only one interconnection point. All models except the DC power flow formulation are limited to radial grids. 

Since the DC model is too inaccurate for distribution grid analysis and the other approaches are limited to one interconnection, we employ a Taylor series linearization of the AC power flow equations using polar coordinates~\cite{Molzahn2019}, which doesn't require the grid to be radial. In the instance used here, there is no representation for the current. However, the formulation can be extended to include current information~\cite{leveringhaus2014}.

Assume that we have a balanced grid and zero line-charging capacitances.
The associated bus-admittance matrix $Y = G+jB \in \mathbb{C}^{|\mathcal N|\times |\mathcal N|}$ is defined as
\[
[Y]_{k,l}=
\begin{cases}
	\sum_{k \in \mathcal N }y_{k,l} & \mbox{if} \quad k=l, \\
	-y_{k,l}, & \mbox{if} \quad k \neq l.
\end{cases}
\]
Here, $y_{k,l} = g_{k,l} + j b_{k,l}\in \mathbb{C}$ is the admittance of  branch $(k,l) \in \mathcal{B}$, and $b_{k,l}$ and $g_{k,l}$ are its susceptance and conductance, respectively.
Note that $y_{k,l} =0$ if $(k,l) \notin \mathcal B$.
The  flow of active power and reactive power  along branch $(k,l) \in \mathcal{B}$ is given by
\begin{subequations} \label{eq:lineFlows_ac}
	\begin{align}
		p_{k,l} &= v_kv_l(G_{k,l}\cos(\theta_{k,l})+B_{k,l}\sin(\theta_{k,l})), \\
		q_{k,l} &= v_k v_l(G_{k,l}\sin(\theta_{k,l})-B_{k,l}\cos(\theta_{k,l})). \label{eq:PFEQq}
	\end{align}
\end{subequations}
Here, $v_k$ is the voltage magnitude at node $k\in \mathcal{N}$ and $\theta_{k,l}\doteq \theta_k - \theta_l$ is the voltage angle difference between the two nodes $k,l\in \mathcal N$.

Using \eqref{eq:lineFlows_ac}, the AC power flow equations for all buses $k \in \mathcal{N}$ read
\begin{align} \label{eq:PFeq}
	&p_k = \sum_{l \in \mathcal{N}} p_{k,l} , \qquad  q_k = \sum_{l \in \mathcal{N}} q_{k,l},
\end{align}
where $p_k, q_k \in \mathbb{R}$ are the net active and reactive power injection at node $k \in \mathcal N$.

The active and reactive power at each bus is defined as
\begin{align} 
    p_k = p_k^g - p_k^d,\\
    q_k = q_k^g -q_k^d \label{eq:nodal_balance}
\end{align} for all $k \in \mathcal N$, where $p_k^g$ and $q_k^g$ represent the active/reactive power feed-in at node $k$ and $p_k^d$ and $q_k^d$ the active/reactive power demand.

To derive a linearized approximation of the AC power flow equations, 
we define $(p_k^\circ, q_k^\circ, v_k^\circ, \theta_k^\circ)$ as the operating point. 
Let $P, Q \in \mathbb{R}^{|\mathcal N|}$ denote the vectors of active and reactive power injections, 
$[p_k]_{k \in \mathcal N}$ and $[q_k]_{k \in \mathcal N}$, respectively, 
and let $V, \theta \in \mathbb{R}^{|\mathcal N|}$ denote the vectors of voltage magnitudes and angles, 
$[v_k]_{k \in \mathcal N}$ and $[\theta_k]_{k \in \mathcal N}$.

Using a first-order Taylor expansion of~\eqref{eq:PFeq} around the operating point $(V^\circ, \theta^\circ)$, 
we obtain
\begin{align} \label{eq:lineFlows}
\begin{bmatrix}
P\\Q
\end{bmatrix} =
\begin{bmatrix}
P^\circ\\Q^\circ
\end{bmatrix} +
J|_{(V^\circ, \theta^\circ)}
\begin{bmatrix}
\theta - \theta^\circ\\[2pt]
V - V^\circ
\end{bmatrix},
\end{align}
where the Jacobian $J$ is given by
\begin{align} \label{eq:jacobian_real}
J =
\begin{bmatrix}
\dfrac{\partial P}{\partial \theta} & \dfrac{\partial P}{\partial V} \\
\dfrac{\partial Q}{\partial \theta} & \dfrac{\partial Q}{\partial V}
\end{bmatrix}.
\end{align}
The submatrices of $J$ represent the sensitivities of active and reactive power injections with respect to 
voltage angles and magnitudes, as defined by the nonlinear AC power flow equations~\eqref{eq:lineFlows_ac}–\eqref{eq:PFeq}. 
Explicit expressions for the entries of $J$ in polar coordinates can be found in~\cite{Molzahn2019, zimmerman2010}. This linearized model provides a convex approximation of the feasible set, used to construct the polyhedral constraint sets $\mathcal X_i$

Besides the power flow constraints, the active power generation must stay within lower and upper bounds. Let $\mathcal G \subseteq \mathcal N$ denote the set of generator buses, then
\begin{align}
	\underline p_k^g \leq p_k^g \leq \bar p_k^g, \qquad k \in \mathcal G. \label{eq:genCstrBds}
\end{align} For ``flexible'' renewable generators, i.e., generators which can alter their set points, we have $\underline p_k^g = 0$ and $\bar p_k = f_k^g$,
where $f_k^g$ represents the maximal power feed in given by irradiation/wind conditions.
Moreover, the renewable generators $k \in \mathcal{G}$ can offer reactive power support $q_k$, constrained by the maximum apparent power $\bar s_k$  and a power factor limit $\alpha$. 
Thus, we consider  affine constraints \cite{contreras2019}
\begin{align} 
	p_k^g \leq \bar s_k\cos(\alpha), \quad -p_k\leq \alpha q_k\leq p_k. \label{eq:ice_cream}
\end{align} 

Additionally, we constrain the voltage magnitude at all buses to stay within bounds
\begin{align} 
	\underline v_k \leq v_k \leq  \bar v_k. \label{eq:voltage_constraint}
\end{align}

Finally, we constrain the active power flowing into the grid at the set of interconnecting branches $(k,l) \in \mathcal B_i^c$, so that both nodes are only importing power from the TSO.
\begin{align}
	p_{k,l} \geq 0. \label{eq:genSlacks}
\end{align} 

The final constraint set for the individual DSO's reads
\begin{align*}
	\mathcal X_i \hspace{-.1cm}\doteq\hspace{-.1cm} \big \{ x_i \in \mathbb{R}^{n_{xi}} \,|\,\eqref{eq:lineFlows}\text{-}\eqref{eq:genSlacks} \text{ hold } \forall k \in \mathcal N_i,\, \text{and } \forall (k,l) \in \mathcal B_i\big \}
\end{align*}
where  $x_i \doteq  \left [ [p_k^g, q_k^g,p_k, q_k, v_k,\theta_k]_{k \in \mathcal N_i},[p_{k,l},q_{k,l}]^\top _{(k,l) \in \mathcal B_i}  ] \right]$ for all $i \in \mathcal C$.
Analogously, $\mathcal X_1$ denotes the feasible set of the TSO, capturing the transmission-level power flow constraints and operational limits.

\subsection{Flexibility aggregation via set projection}
We can now define the coupling variable $z_i^\top \doteq \left [[p_{k,l}, q_{k,l}]_{(k,l) \in \mathcal B_i^c  }, {[v_k,\theta_k]_{k\in \mathcal N_i^c}}\right ]$. The importance of considering the voltage magnitude as a coupling variable has been demonstrated in \cite{Bandeira2024}. There, the voltage angle at the interconnection is considered equal to zero and not a coupling variable. For the AC power flow equations, the active and reactive power on branch $(k,l) \in  \mathcal B$ depends only on the voltage angle difference $(\theta_k-\theta_l)$. We have $p_{k,l}(v_k,v_l, \theta_k,\theta_l) = p_{k,l}(v_k,v_l, \theta_k + \bar \theta,\theta_l+\bar \theta)$ for any $\bar \theta \in \mathbb R$, so setting $\theta_k = 0$ in $ \mathcal X_i$  does not alter the solution space. Here, no such simplification can be made. We must consider the angle difference between the interconnections, as shown in~\cite{stark2024}, but not their actual value. The difference is invariant w.r.t. a constant offset.

We propose to look at aggregation as equivalent to a set projection of the feasible set of the DSOs, $\mathcal X_i$, on the coupling variables $z_i$.
The \emph{set projection} of $\mathcal X \subseteq \mathcal Y \times \mathcal Z$ onto $\mathcal Z$ is defined as~\cite{Rakovic2006}
 \begin{align}\label{eq:setProj}
 	\operatorname{proj}_{\mathcal Z}(\mathcal X) \doteq \{ z \in \mathcal Z  \;|\;\exists \; y \in \mathcal Y \text{ with  } (z,y) \in \mathcal X\} \subseteq \mathcal Z.
 \end{align}
The above definition \eqref{eq:setProj} can be interpreted as searching for all coupling variables to the TSO, $z_i$, for which there exist local variables $y_i$, i.e., admissible set-points of the local DERs, which satisfy all grid constraints.
 Hence, aggregation can be interpreted as computing the FOR
 \begin{align} \label{eq:proj}
 	\mathcal P_i \doteq \operatorname{proj}_{ \mathcal Z_i}{ {\mathcal  X}}_i.
 \end{align}

There are many advantages to computing the FOR using set projection. First, it is a flexible approach, since we can define the coupling variable depending on the task at hand, without altering the methodology or the framework. Additionally, if using polytope set projection, the resulting FOR can be expressed as the polytope $A_i z_i \leq b_i$ . The DSO then only needs to communicate the matrices $A_i,b_i$ to the TSO, which can then easily add the linear inequalities to its own optimization problem, without adding much complexity.
Here, we focus on the flexibility aggregation of one distribution grid with multiple interconnections, i.e. the computation of the FOR $\mathcal P_i$.

\section{Simulation} \label{sec:simulation}
\begin{figure}[t]
	\centering
	\includegraphics[angle=-90, width=0.25\textwidth,trim={0 0 12cm 0},clip]{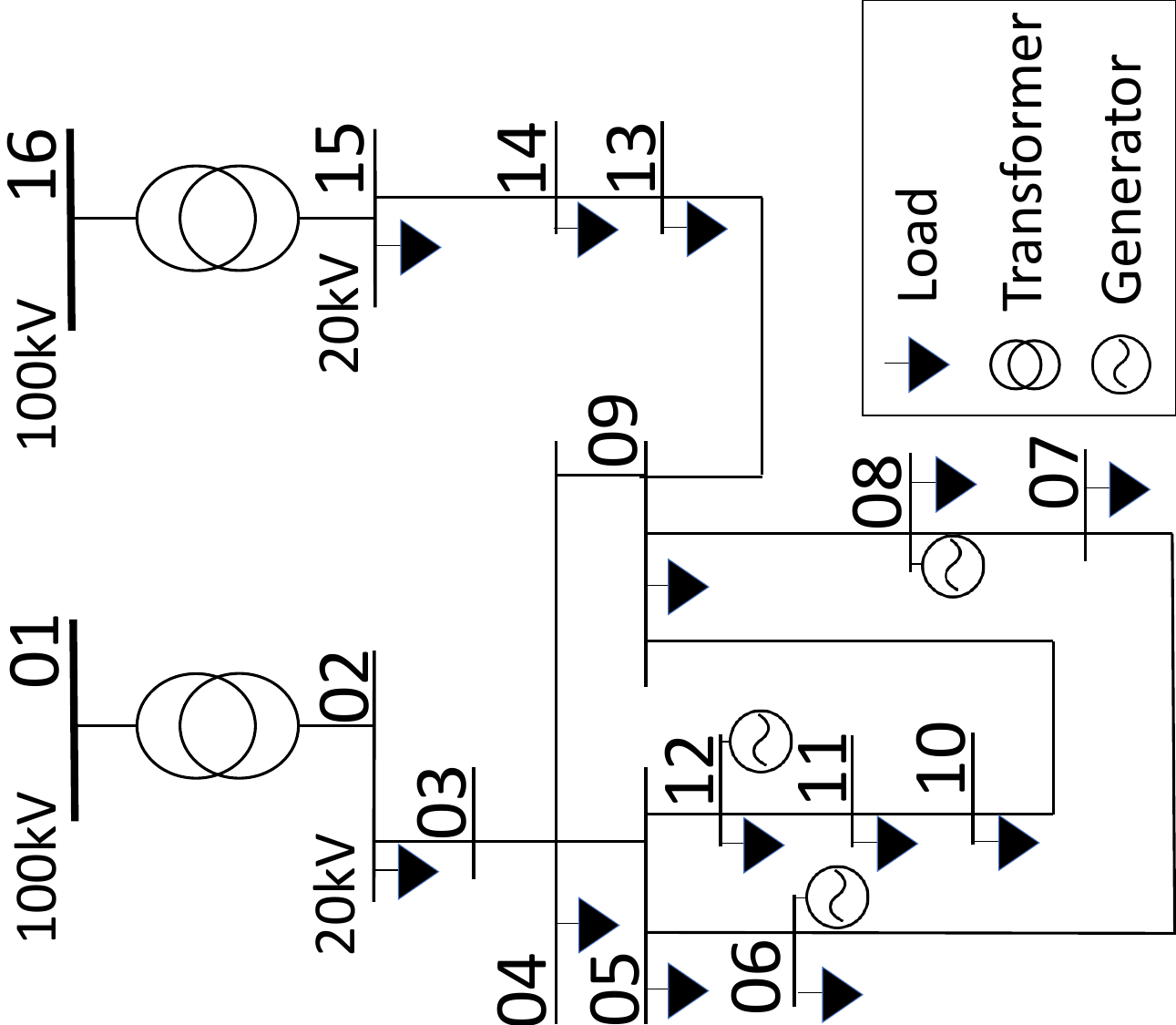}
	\caption{CIGRE medium voltage European grid \cite{cigre}.}
	\label{fig:grid}
\end{figure}

For the simulation, we employ the CIGRE Medium Voltage European Grid testbench~\cite{cigre}, which represents a meshed distribution network with two transformers connecting to the high-voltage grid, as illustrated in Figure~\ref{fig:grid}. Renewable generators are placed at nodes 6, 8, and 12, and their active power feed-in capacities are doubled relative to the original specification. Following the setup in the original benchmark, we first consider the case in which both interconnection points connect to the same high-voltage bus, i.e., node~1 = node~16.

\subsection{Suitability of approximation}
Firstly, we evaluate the suitability of the linearized formulation introduced in Section~\ref{sec:problem} for the purpose of flexibility aggregation in the meshed distribution grid. Specifically, we compare the resulting FOR with that obtained using an optimization-based sampling method based on the nonlinear AC power flow formulation introduced in the literature.

To this end, we compute samples $z_k^s = (p_{k,l}, q_{k,l})$ by solving multiple optimization problems of the form
\begin{align} \label{eq:boundaryProb}
z_k^t = \arg \min_{z_i, y_i}  c_l^\top z_k
\quad \text{subject to} \quad (y_i, z_i) \in \mathcal X_i^{AC},
\end{align}
where we assume the voltage magnitude and angle reference
\begin{align}
v_k = 1\,\text{p.u.}, \quad \theta_k = 0,
\label{eq:reference}
\end{align}
and define
\begin{align*} \label{eq:setAC}
\mathcal X_i^{AC} \doteq
\Big\{ 
x_i \in \mathbb{R}^{n_{x_i}} \;&\Big|\;
\eqref{eq:lineFlows_ac}\text{-}\eqref{eq:nodal_balance},\;
\eqref{eq:genCstrBds}\text{-}\eqref{eq:genSlacks},\;
\text{and}~\eqref{eq:reference}\\ \qquad &\quad\text{hold } 
\forall k \in \mathcal N_i, \; (k,l) \in \mathcal B_i
\Big\}.
\end{align*}

Different boundary points of the FOR are obtained by varying the cost vector $c_l \in {-1,0,1}^{n_z}$
to explore the extreme points~\cite{capitanescu2018,silva2018}. To refine the sampling, the resulting region is discretized, and one variable component is fixed for each grid point before re-solving~\eqref{eq:boundaryProb}. A detailed description of this process is provided in~\cite{Engelmann2025}. The resulting FOR is defined as
\begin{align*}
\mathcal F_i = \operatorname{conv}(\{z_k^t\}_{t\in \mathcal T}),
\end{align*}
where $\mathcal{T}$ denotes the set of all sample indices and $\operatorname{conv}(\mathcal W)$ denotes the convex hull over a set of points $\mathcal W$.

\begin{figure}[t]
	\centering
	\includegraphics[width=0.95\linewidth]{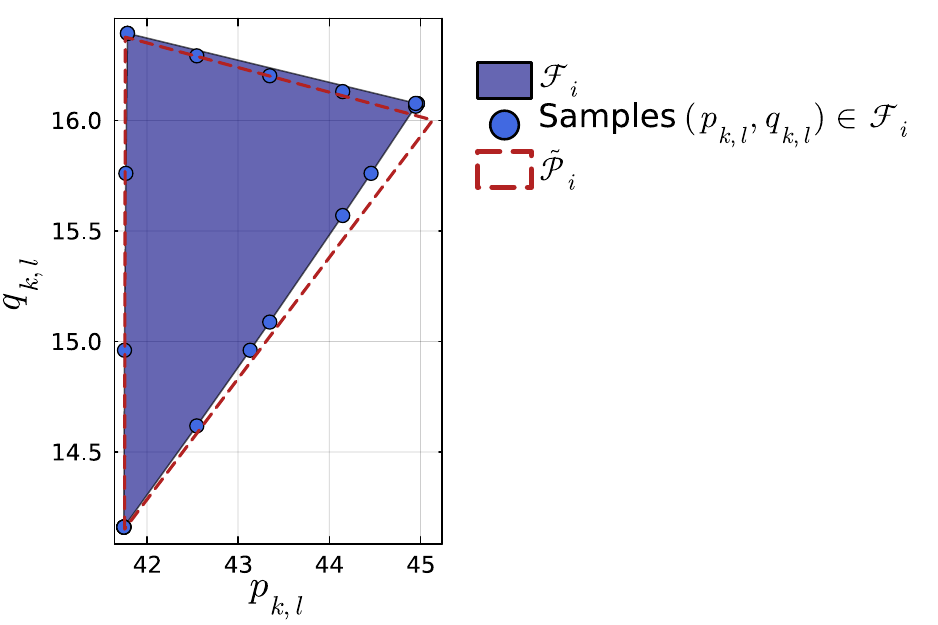}
	\caption{Comparison of FORs obtained with AC-OPF sampling and the proposed linearized set projection method for a single interconnection point.}
	\label{fig:1slack}
\end{figure}

Subsequently, we compute the linearized FOR $\tilde{\mathcal P}i$ as described in Section~\ref{sec:problem}, 
\begin{align*}
\tilde{\mathcal P}_i = \operatorname{proj}_{z_i}(\tilde{\mathcal X}_i)
\end{align*} 
with $z_i = [p_{k,l}, q_{k,l}]^\top$ and the extended constraint set
\begin{align*} \tilde{\mathcal X}_i \doteq \Big\{ x_i \in \mathcal X_i \;\big|\; \text{\eqref{eq:reference} holds} \Big\}.\end{align*} The operating point is chosen as the solution of \begin{align} 
(x_i^\circ, z_i^\circ) = &\arg \min_{x_i, z_i} \;
\sum_{k \in \mathcal{G}} \Bigl[ -c_1 (p_k^g)^2 + c_2 (q_k^g)^2 \Bigr] \\
& \text{subject to} \quad 
(x_i, z_i) \in \mathcal{X}_i^{AC}.
\end{align} corresponding to the case that minimizes the power curtailment of renewable generators. The resulting regions are depicted in Figure~\ref{fig:1slack}. The linearized formulation approximates the FOR obtained with the AC model closely, although it includes certain infeasible points that may lead to infeasibility in the coupled optimization problem~\cite{Engelmann2025}. 

\subsection{Flexibility aggregation via set projection}

Now, with a reasonably good approximation for the one interconnection bus case, we portray the interconnection buses as in Figure~\ref{fig:grid}, as two independent buses, $k = 1$ and $k = 16$, each connected to one of the different transformers. We then extend the coupling variables for the multi-interconnections case as the 8-dimensional vector $z_i \doteq  \left [ p_{1,2}, q_{1,2}, v_1, \theta_1, p_{16,15}, q_{16,15}, v_{16}, \theta_{16}  \right]$. Since in this case we have only two interconnections, and as explained in Section \ref{sec:problem} the difference between the angles is what matters, we set $\theta_0 = 0$ and use $\theta_1-\theta_{16}$ as the coupling variable, reducing the coupling vector to 7 dimensions.

By projecting the full constraint set $\mathcal X_i$ into the coupling variable $z_i$, we obtain a FOR, which is a 7-dimensional polyhedron, representing the flexibility available at each interconnection as well as the dependent relationship between the variables.

\subsection{Results discussion}
To interpret the abstract results in a more tangible way, we illustrate the relationship between pairs of variables, using the set projection, of the 7D-polyhedron into the desired two variables. 

\begin{figure}[htbp]
	\centering
	\begin{subfigure}[b]{0.5\textwidth}
		\centering
	\includegraphics[width=\textwidth]{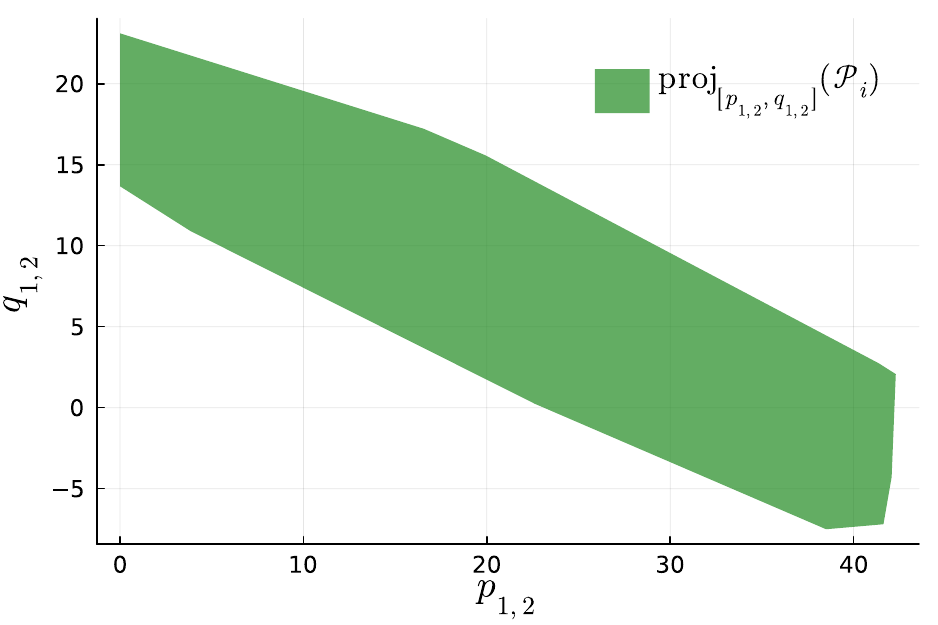}
		\caption{Active and reactive power at interconnection on line $(1,2)$ }
		\label{fig:pq1}
	\end{subfigure}
	\hfill
	\begin{subfigure}[b]{0.5\textwidth}
		\centering
		\includegraphics[width=\textwidth]{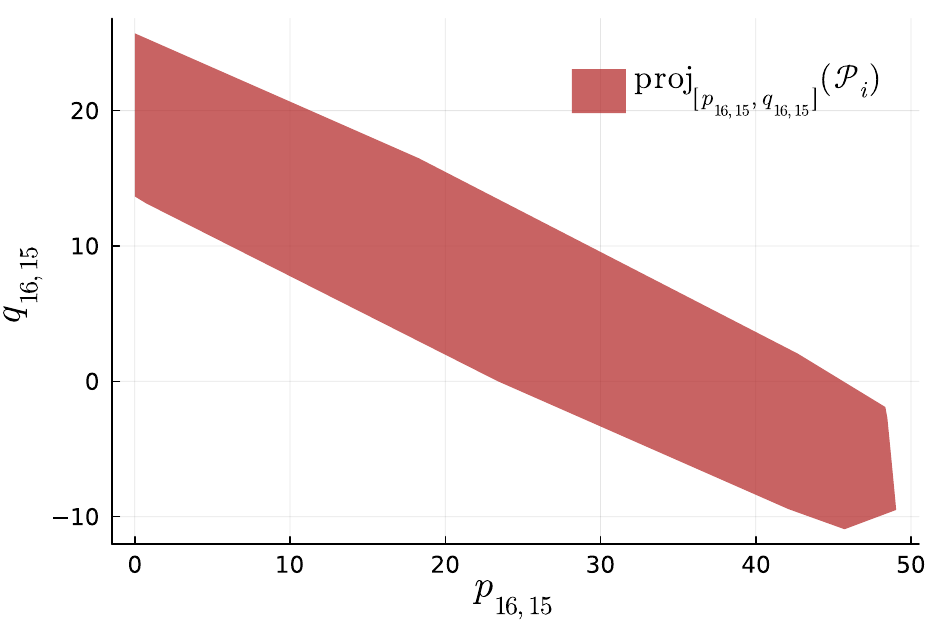}
		\caption{Active and reactive power at interconnection on line $(16,15)$ }
		\label{fig:pq2}
	\end{subfigure}
	\caption{Projected views for active and reactive power of the 7D-FOR for the case with two interconnection points.}
	\label{fig:pq}
\end{figure}

First, we analyze the aggregated active and reactive power at the individual interconnections. As shown in Figure~\ref{fig:pq}, the available maximum power at each interconnection differs due to the location of the assets. It is important to note that these results should not be considered in isolation—neither relative to each other nor to the remaining variables—as all are interconnected. In other words, these charts cannot be interpreted independently.

\begin{figure}[t]
	\centering
	\includegraphics[width=0.95\linewidth]{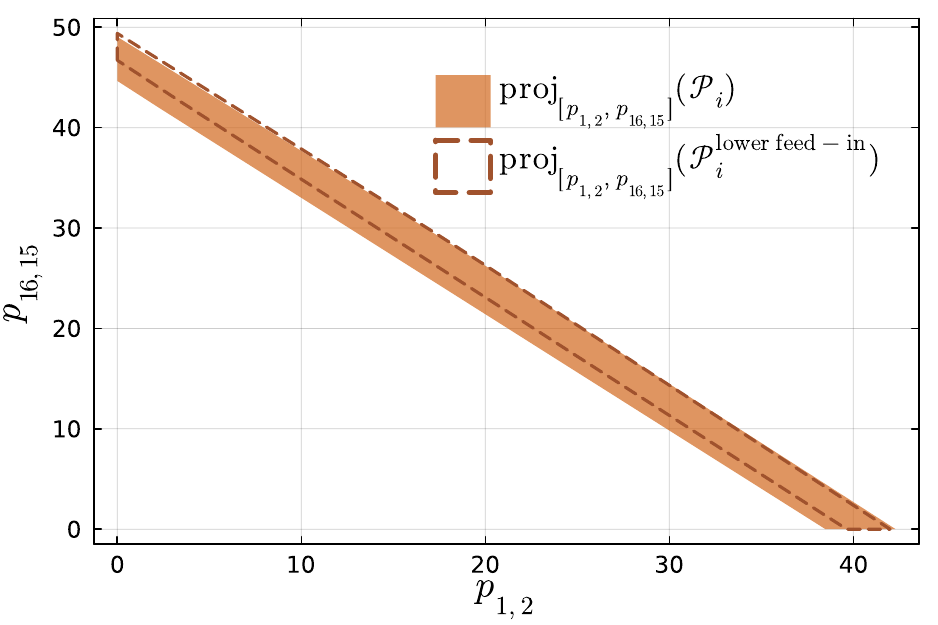}
	\caption{Visualization of the relationship between active power flows  $p_{1,2}$ and  $p_{16,15}$ based on the 7D-FOR projection.}
	\label{fig:pp}
\end{figure}

To better illustrate this point, we analyze the dependency of the active power in one interconnection on the active power in the second interconnection. More specifically, Figure~\ref{fig:pp} illustrates a negative dependency between the two active powers: as the active power flowing into the distribution grid increases at one interconnection, it decreases at the other. 
This indicates that the power demand can be met through both interconnections, such that a higher export on one side reduces the power supplied from the other.
The interior of the FOR reflects the inherent flexibilities within the distribution grid. The dashed line represents the area where the renewable generators' power feed-in is not increased, illustrating the reduced flexibility within the grid.

Next, we compare the active and reactive power requirements to the single-interconnection case, as shown in Figure~\ref{fig:1slack}. To this end, we define the aggregated variables
\begin{align*}
    p_{\text{sum}} = p_{1,2} + p_{16,15}, \quad  
    q_{\text{sum}} = q_{1,2} + q_{16,15},
\end{align*}
and compute the corresponding FOR via projection as
\begin{align*}
    \mathcal P_i' = \operatorname{proj}_{[p_{\text{sum}},q_{\text{sum}},v_1,v_{16},\theta_1-\theta_{16}]} (\mathcal X_i),
\end{align*}
where $\mathcal X_i$ is the extended constraint set including the voltage and power flow constraints at both interconnection buses.

Since the grid topology now includes two independent interconnections, the FOR depends not only on the total power exchange but also on the individual voltage magnitudes $v_1$, $v_{16}$ and the relative voltage angle $\theta_1-\theta_{16}$. Therefore, we consider three distinct case setups to isolate the effects of these additional degrees of freedom:

\begin{table}[h!]
\centering
\caption{Overview of the considered cases and corresponding constraints on the interconnection buses.}
\label{tab:cases}
\begin{tabular}{c c c}
\toprule
\textbf{Case} & \textbf{Voltage magnitudes} & \textbf{Angle difference} \\
\midrule
1 & $v_1$, $v_{16}$ free & $\theta_1 - \theta_{16}$ free \\
2 & $v_1$, $v_{16}$ free & $\theta_1 - \theta_{16} = 0$ \\
3 & $v_1 = v_{16} = 1\,\text{p.u.}$ & $\theta_1 - \theta_{16} = 0$ \\
\bottomrule
\end{tabular}
\end{table}

\begin{figure}[t]
	\centering
	\includegraphics[width=0.95\linewidth]{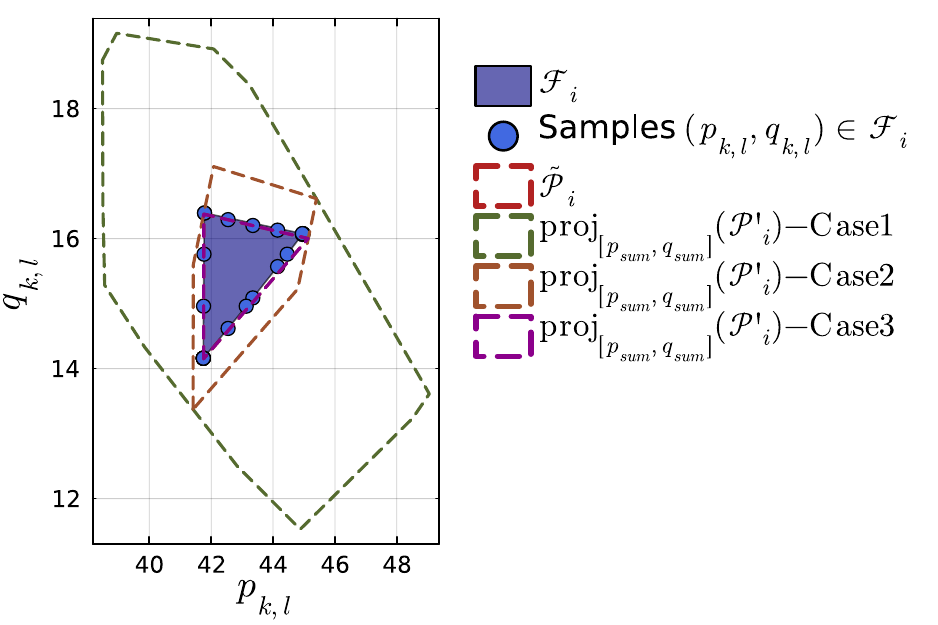}
	\caption{Total active and reactive power requirements for the three cases of the two-interconnection scenario, compared to the single-interconnection case.}
	\label{fig:theta_taylor}
\end{figure}

In Case~1, all interconnection variables are free, resulting in the largest FOR. Case~2 constrains the angle difference, effectively reducing the dimensionality of the feasible region and leading to a smaller FOR, shown in Figure~\ref{fig:theta_taylor} in orange. Finally, Case~3 fixes both voltage magnitudes and the angle difference, producing a region comparable to the single-interconnection scenario (Figure~\ref{fig:1slack}). These comparisons illustrate how the additional degrees of freedom in multi-interconnection grids influence the aggregated flexibility and highlight the importance of including both voltage and angle variations in the FOR computation.

These observations are directly explainable from a power systems perspective. In Case~3, where the voltage angles and magnitudes are fixed, we essentially recover the scenario in which both interconnections are connected to the same bus, yielding a FOR identical to that in Figure~\ref{fig:1slack}. In the other cases, the voltage angles and magnitudes are allowed to vary, and in the extreme case where they are completely free (Case~1).

From a TSO perspective, this behavior is intuitive: multiple interconnections enable the distribution grid to be used for both flexibility provision and power transfer between buses. It is important to note, however, that in a real transmission system, the voltage angles and magnitudes at the interconnection buses are not entirely free; they are influenced by the connected TSO grid, its electrical parameters, and the setpoints of the deployed equipment.

\begin{figure}[t]
	\centering
	\includegraphics[width=0.95\linewidth]{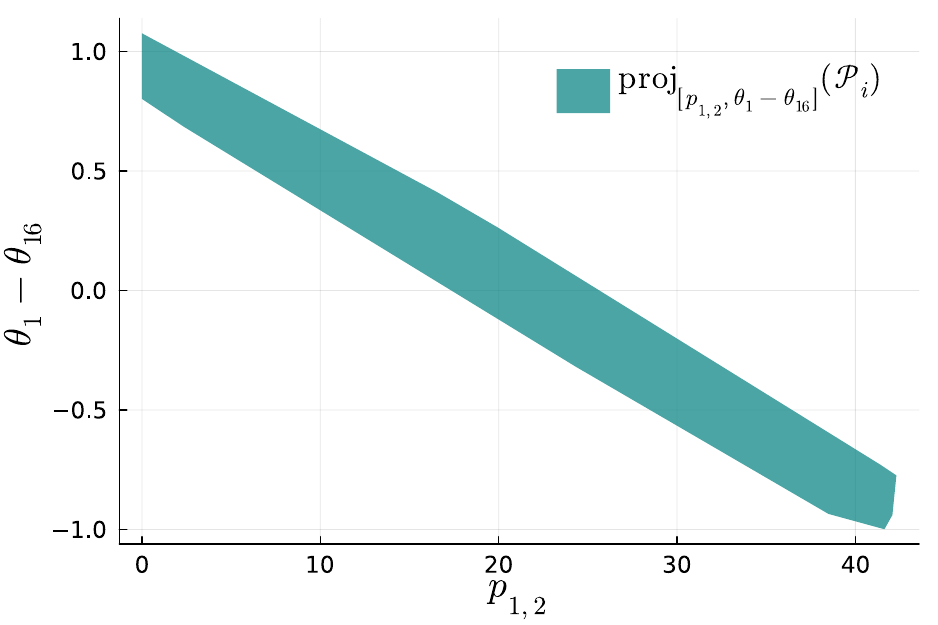}
	\caption{Visualization of the relationship between active power $p_{1,2}$ and the voltage angle difference $\theta_1 - \theta_{16}$ based on the 7D-FOR projection.}
	\label{fig:p1_difftheta}
\end{figure}

The dependency between the active power at one interconnection and the voltage angle difference can be visualized by projecting our 7-dimensional FOR onto the variables $p_{1,2}$ and $\theta_1 - \theta_{16}$, as shown in Figure~\ref{fig:p1_difftheta}. The projection highlights how variations in the angle difference directly influence the feasible range of active power that can be drawn from the distribution grid at this interconnection. A similar relationship is observed for the other interconnection, $p_{16,15}$, reflecting the inherent coupling between the interconnections in the network. These projections provide an intuitive view of how the flexibility of the distribution grid is constrained by both network topology and voltage conditions.

This dependency may also help explain the differences observed in the computed regions in \cite{silva2023,stark2024}. In that work, the flexibility areas obtained using an equivalent circuit for the transmission or distribution grid are smaller than those computed for the distribution grid in isolation. This can be attributed to the natural restriction of the voltage differences imposed by the physical properties of the transmission grid. In contrast, our aggregated 7D-FOR explicitly accounts for the voltage magnitude and angle coupling. As a result, when the TSO receives the 7D-FOR, it will inherently operate within the subset of the region that is feasible according to its own grid constraints, making an equivalent circuit representation unnecessary.

\begin{figure}[t]
	\centering
	\includegraphics[width=0.95\linewidth]{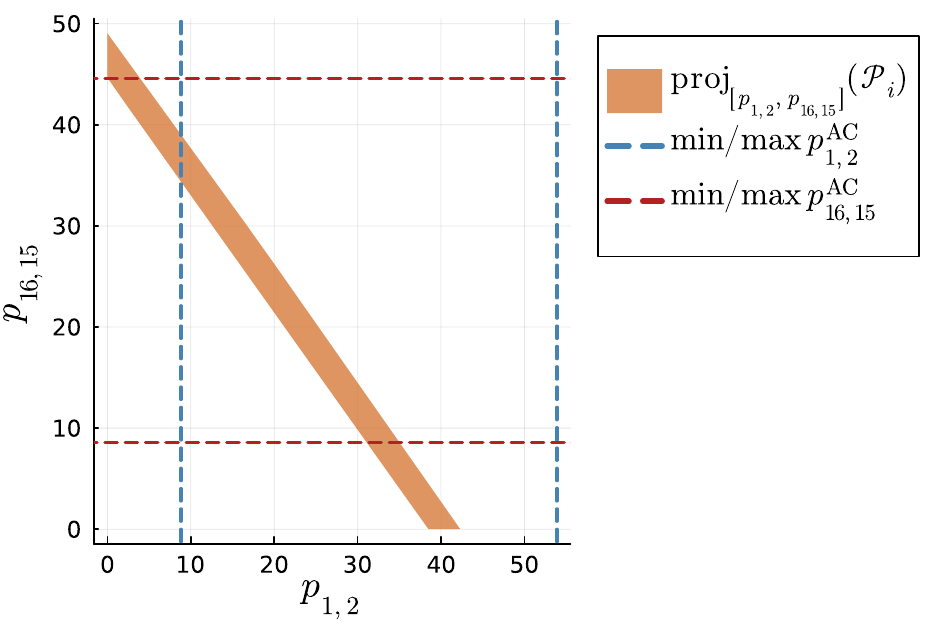}
	\caption{Projection of the 7D-FOR onto active powers $p_{1,2}$ and $p_{16,15}$ for the two-interconnection scenario, compared to AC-OPF extreme values sampling results for the same case.}
	\label{fig:comp_ac}
\end{figure}

Lastly, we examine the potential modeling inaccuracies of the chosen power flow formulation due to linearization and its implications on the high-dimensional FOR. Although we cannot easily compute the 7D-FOR with the AC power flow, we can solve the optimization problem for the maximum values to get an idea of the difference. 
By treating multiple variables as coupling variables, the resulting operating points may lie farther from the linearization points, which inevitably increases potential linearization errors near the boundary of the high-dimensional FOR, as illustrated in Figure~\ref{fig:comp_ac}. The actual impact of these errors on system operation, however, depends on the TSO’s decisions: points near the FOR boundary may introduce inaccuracies, whereas points closer to the linearization point are expected to be less critical.

\section{Conclusion and Outlook} \label{sec:conclusion}
We have presented a flexibility aggregation scheme for distribution grids with multiple interconnections. We have extended works based on set projections to be applied to grids with two or more interconnections. The resulting multidimensional FOR helps to gain insight into the relationship between the physical variables at both interconnections and can be easily communicated to the upper grid as a polytope.\\
In future work, we will look into different linear formulations for meshed grids that can be used for set projection and model the behavior of the grid outside of the operational point better.\\
A further aspect is to make the projection step itself more efficient, e.g., by using inner approximations of the FOR, which are easier to project.

\renewcommand*{\bibfont}{\small}
\printbibliography

@article{Engelmann2025,
  title = {Approximate {{Dynamic Programming With Feasibility Guarantees}}},
  author = {Engelmann, Alexander and Bandeira, Ma{\'i}sa Beraldo and Faulwasser, Timm},
  year = {2025},
  journal = {IEEE Transactions on Control of Network Systems},
  pages = {1565 -- 1576},
  issn = {2325-5870},
  doi = {10.1109/TCNS.2025.3526715},
  urldate = {2025-01-15}
}

@article{capitanescu2018,
  title = {{{TSO}}\textendash{{DSO}} Interaction: {{Active}} Distribution Network Power Chart for {{TSO}} Ancillary Services Provision},
  shorttitle = {{{TSO}}\textendash{{DSO}} Interaction},
  author = {Capitanescu, Florin},
  year = {2018},
  month = oct,
  journal = {Electric Power Systems Research},
  volume = {163},
  pages = {226--230},
  issn = {0378-7796},
  doi = {10.1016/j.epsr.2018.06.009},
  langid = {english}
}

@inproceedings{mayorga2018,
  title = {Determination of the {{Time-Dependent Flexibility}} of {{Active Distribution Networks}} to {{Control Their TSO-DSO Interconnection Power Flow}}},
  booktitle = {{{Power Systems Computation Conference}} ({{PSCC}})},
  author = {Mayorga Gonzalez, D. and Hachenberger, J. and Hinker, J. and Rewald, F. and H{\"a}ger, U. and Rehtanz, C. and Myrzik, J.},
  year = {2018},
  month = jun,
  pages = {1--8},
  doi = {10.23919/PSCC.2018.8442865}
}

@article{silva2018,
  title = {Estimating the Active and Reactive Power Flexibility Area at the {{TSO-DSO}} Interface},
  author = {Silva, Jo{\~a}o and Sumaili, Jean and Bessa, Ricardo J. and Seca, Lu{\'i}s and Matos, Manuel A. and Miranda, Vladimiro and Caujolle, Mathieu and Goncer, Bel{\'e}n and {Sebastian-Viana}, Maria},
  year = {2018},
  month = sep,
  journal = {IEEE Transactions on Power Systems},
  volume = {33},
  number = {5},
  pages = {4741--4750},
  issn = {1558-0679},
  doi = {10.1109/TPWRS.2018.2805765}
}

@inproceedings{contreras2019,
  title = {Time-{{Based Aggregation}} of {{Flexibility}} at the {{TSO-DSO Interconnection Point}}},
  booktitle = {{{IEEE Power}} \& {{Energy Society General Meeting}} ({{PESGM}})},
  author = {Contreras, Daniel A. and Rudion, Krzysztof},
  year = {2019},
  month = aug,
  pages = {1--5},
  issn = {1944-9933},
  doi = {10.1109/PESGM40551.2019.8973421}
}

@article{sarstedt2021,
  title = {Survey and {{Comparison}} of {{Optimization-Based Aggregation Methods}} for the {{Determination}} of the {{Flexibility Potentials}} at {{Vertical System Interconnections}}},
  author = {Sarstedt, Marcel and Klu{\ss}, Leonard and Gerster, Johannes and Meldau, Tobias and Hofmann, Lutz},
  year = {2021},
  month = jan,
  journal = {Energies},
  volume = {14},
  number = {3},
  pages = {687},
  publisher = {{Multidisciplinary Digital Publishing Institute}},
  issn = {1996-1073},
  doi = {10.3390/en14030687},
  urldate = {2023-07-26},
  copyright = {http://creativecommons.org/licenses/by/3.0/},
  langid = {english}
}

@misc{majumdar2023,
  title = {Distribution Grid Power Flexibility Aggregation at Multiple Interconnections between the High and Extra High Voltage Grid Levels},
  author = {Majumdar, Neelotpal and Sarstedt, Marcel and Hofmann, Lutz},
  year = {2023},
  month = mar,
  number = {arXiv:2303.01107},
  eprint = {2303.01107},
  %primaryclass = {cs, eess},
  publisher = {{arXiv}},
  doi = {10.48550/arXiv.2303.01107},
  urldate = {2023-09-30},
  abstract = {The energy transition towards renewable based power provision requires improved monitoring and control of distributed energy resources (DERs), installed predominantly at the distribution grid level. Due to the gradual phase out of thermal generation, a shift of ancillary services provision like voltage control, congestion management and dynamic support from DERs is underway. Increased planning for procurement of ancillary services from underlying grid levels is required. Therefore, provision of flexible active and reactive power potentials from distribution system operators to transmission system operators at the vertical system interface is a subject of current research. At present, provision of active and reactive power flexibilities (PQ-flexibilities) across radial system interconnections are investigated, which involves a single transformer branch interconnection across two different voltage levels. Inclusion of multiple interconnections in a meshed grid structure increases the complexity as proximal interdependencies of interconnections to PQ-flexibilities require consideration. The objective of this paper is to address the flexibility aggregation across multiple vertical interconnections. Alternating current power transfer distribution factors (AC-PTDFs) are used to determine the power flow across the interconnections. Subsequent integration in a linear optimization environment controls the interconnection power flows (IPF) using a weighted objective function. Therefore, power flow regulation is enabled according to the requirements and specifications of both the underlying and overlaying grid level. The results show interdependent concentric flexibility regions or Feasible Operating Regions (FORs) in accordance with the manipulation of the weighted objective function.},
  archiveprefix = {arxiv},
  keywords = {Electrical Engineering and Systems Science - Systems and Control},
  file = {C\:\\Users\\smmsbera\\Zotero\\storage\\QTA9K4C4\\Majumdar et al. - 2023 - Distribution grid power flexibility aggregation at.pdf;C\:\\Users\\smmsbera\\Zotero\\storage\\DRJ3J4AJ\\2303.html}
}

@inproceedings{sarstedt2020,
  title = {Simulation of {{Hierarchical Multi-Level Grid Control Strategies}}},
  booktitle = {2020 {{International Conference}} on {{Smart Grids}} and {{Energy Systems}} ({{SGES}})},
  author = {Sarstedt, Marcel and Klu\ss{}, Leonard and Dokus, Marc and Hofmann, Lutz and Gerster, Johannes},
  year = {2020},
  month = nov,
  pages = {175--180},
  publisher = {{IEEE}},
  address = {{Perth, Australia}},
  doi = {10.1109/SGES51519.2020.00038},
  urldate = {2023-09-30},
  isbn = {978-1-72818-550-7},
  langid = {english}
}

@article{Bandeira2024,
  title = {An {{ADP Framework}} for {{Flexibility}} and {{Cost Aggregation}}: {{Guarantees}} and {{Open Problems}}},
author = {Beraldo Bandeira, Ma{\'i}sa and Faulwasser, Timm and Engelmann, Alexander},
  year = {2024},
  month = sep,
volume = {234},
  pages = {110818},
  journal = {Electric Power Systems Research},
  issn = {0378-7796},
  doi = {10.1016/j.epsr.2024.110818},
  urldate = {2025-01-15}
}

@misc{Bandeira2025,
  title = {Complexity {{Reduction}} for {{TSO-DSO Coordination}}: {{Flexibility Aggregation}} vs. {{Distributed Optimization}}},
  shorttitle = {Complexity {{Reduction}} for {{TSO-DSO Coordination}}},
  author = {Bandeira, Ma{\'i}sa Beraldo and Engelmann, Alexander and Faulwasser, Timm},
  year = {2025},
  month = sep,
  number = {arXiv:2509.10595},
  eprint = {2509.10595},
  primaryclass = {eess},
  publisher = {arXiv},
  doi = {10.48550/arXiv.2509.10595},
  urldate = {2025-09-16},
  archiveprefix = {arXiv},
}

@article{silva2023,
  title = {A {{Data-Driven Approach}} to {{Estimate}} the {{Flexibility Maps}} in {{Multiple TSO-DSO Connections}}},
  author = {Silva, Jo{\~a}o and Sumaili, Jean and Silva, Bernardo and Carvalho, Leonel and Retorta, F{\'a}bio and Staudt, Maik and Miranda, Vladimiro},
  year = {2023},
  month = mar,
  journal = {IEEE Transactions on Power Systems},
  volume = {38},
  number = {2},
  pages = {1908--1919},
  issn = {1558-0679},
  doi = {10.1109/TPWRS.2022.3214106},
  urldate = {2024-02-15}
}

@book{cigre,
  title = {Benchmark Systems for Network Integration of Renewable and Distributed Energy Resources},
  editor = {{Conseil international des grands r{\'e}seaux {\'e}lectriques}},
  year = {2014},
  publisher = {CIGR{\'E}},
  address = {Paris},
  isbn = {978-2-85873-270-8},
  langid = {english},
  lccn = {620}
}

@techreport{zimmerman2010,
  title = {{{AC Power Flows}}, {{Generalized OPF Costs}} and Their {{Derivatives}} Using {{Complex Matrix Notation}}},
  author = {Zimmerman, Ray D},
  year = {2010},
  langid = {english},
    institution = {Matpower}
}

@inproceedings{stark2024,
  title = {Novel {{Approach}} for {{Flexibility Aggregation}} at {{Multiple Vertical TSO-DSO Interconnections}}},
  booktitle = {Transformation Der {{Stromversorgung}} -- {{Netzregelung}} Und {{Systemf{\"u}hrung}}; 15. {{ETG}}/{{GMA-Fachtagung}}, {{Netzregelung}} Und {{Systemf{\"u}hrung}}``},
  author = {Stark, Lars and Hofmann, Lutz},
  year = {2024},
  month = sep,
  pages = {18--24},
  urldate = {2025-10-13}
}

@inproceedings{lyu2023,
  title = {Feasible {{Set Projection Method}} for {{Calculating Multi-Line-Coupled ATC}}},
  booktitle = {2023 3rd {{Power System}} and {{Green Energy Conference}} ({{PSGEC}})},
  author = {Lyu, Chengming and He, Letian and Tan, Zhenfei and Lin, Ye and Luo, Xi and Yan, Zheng},
  year = {2023},
  month = aug,
  pages = {791--797},
  doi = {10.1109/PSGEC58411.2023.10255986},
  urldate = {2025-10-13}
}

@article{tan2024,
  title = {Non-{{Iterative Solution}} for {{Coordinated Optimal Dispatch}} via {{Equivalent Projection}}---{{Part I}}: {{Theory}}},
  shorttitle = {Non-{{Iterative Solution}} for {{Coordinated Optimal Dispatch}} via {{Equivalent Projection}}---{{Part I}}},
  author = {Tan, Zhenfei and Yan, Zheng and Zhong, Haiwang and Xia, Qing},
  year = {2024},
  month = jan,
  journal = {IEEE Transactions on Power Systems},
  volume = {39},
  number = {1},
  pages = {890--898},
  issn = {1558-0679},
  doi = {10.1109/TPWRS.2023.3258066},
  urldate = {2025-10-13}
}

@inproceedings{wang2023,
  title = {A {{Projection-Based Approach}} for {{Distributed Energy Resources Aggregation}}},
  booktitle = {2023 {{IEEE PES Innovative Smart Grid Technologies Europe}} ({{ISGT EUROPE}})},
  author = {Wang, Yiran and Zhong, Haiwang and Ruan, Guangchun},
  year = {2023},
  month = oct,
  pages = {1--5},
  doi = {10.1109/ISGTEUROPE56780.2023.10408587},
  urldate = {2025-10-13}
}

@inproceedings{leveringhaus2014,
  title = {Comparison of Methods for State Prediction: {{Power Flow Decomposition}} ({{PFD}}), {{AC Power Transfer Distribution}} Factors ({{AC-PTDFs}}), and {{Power Transfer Distribution}} Factors ({{PTDFs}})},
  shorttitle = {Comparison of Methods for State Prediction},
  booktitle = {2014 {{IEEE PES Asia-Pacific Power}} and {{Energy Engineering Conference}} ({{APPEEC}})},
  author = {Leveringhaus, T. and Hofmann, L.},
  year = {2014},
  month = dec,
  pages = {1--6},
  issn = {2157-4847},
  doi = {10.1109/APPEEC.2014.7066183},
  urldate = {2025-10-13}
}

@inproceedings{stark2023,
  title = {Determination of {{Interdependent Feasible Operation Regions}} at {{Multiple TSO-DSO Interconnections}}},
  booktitle = {2023 {{IEEE PES Innovative Smart Grid Technologies Europe}} ({{ISGT EUROPE}})},
  author = {Stark, Lars and Sarstedt, Marcel and Hofmann, Lutz},
  year = 2023,
  month = oct,
  pages = {1--6},
  doi = {10.1109/ISGTEUROPE56780.2023.10408712},
  urldate = {2024-02-15}
}

@inproceedings{MayorgaGonzalez2018,
  title = {Determination of the {{Time-Dependent Flexibility}} of {{Active Distribution Networks}} to {{Control Their TSO-DSO Interconnection Power Flow}}},
  booktitle = {2018 {{Power Systems Computation Conference}} ({{PSCC}})},
  author = {Mayorga Gonzalez, D. and Hachenberger, J. and Hinker, J. and Rewald, F. and H{\"a}ger, U. and Rehtanz, C. and Myrzik, J.},
  year = {2018},
  month = jun,
  pages = {1--8}
}

@article{Molzahn2019,
  title = {A {{Survey}} of {{Relaxations}} and {{Approximations}} of the {{Power Flow Equations}}},
  author = {Molzahn, Daniel K. and Hiskens, Ian A.},
  year = {2019},
  month = feb,
  journal = {Foundations and Trends{\textregistered} in Electric Energy Systems},
  volume = {4},
  number = {1-2},
  pages = {1--221},
  publisher = {Now Publishers, Inc.},
  langid = {english}
}

@article{Rakovic2006,
  title = {Reachability Analysis of Discrete-Time Systems with Disturbances},
  author = {Rakovic, S.V. and Kerrigan, E.C. and Mayne, D.Q. and Lygeros, J.},
  year = {2006},
  month = apr,
  journal = {IEEE Transactions on Automatic Control},
  volume = {51},
  number = {4},
  pages = {546--561}
}

\end{document}